# PROJECTIVE Q-FACTORIAL TORIC VARIETIES COVERED BY LINES

CINZIA CASAGRANDE AND SANDRA DI ROCCO

ABSTRACT. We give a structural theorem for Q-factorial toric varieties covered by lines in $\mathbb{P}^N$, and compute their dual defect. This yields a characterization of defective Q-factorial toric varieties in $\mathbb{P}^N$. The combinatorial description of such varieties is used to characterize some finite sets of monomials with discriminant equal to one.

## INTRODUCTION

Let $X \subset \mathbb{P}^N$ be a projective variety such that through any point of $X$ there is a line contained in $X$. Such $X$ is said to be covered by lines, or uniruled by lines. The main result of this paper is a structural theorem for Q-factorial toric varieties covered by lines.

**Theorem 4.2.** *Let $X \subset \mathbb{P}^N$ be a Q-factorial toric variety, covered by lines. Then there exist (and are uniquely determined) a projective, Q-factorial toric variety $Z$, and a surjective and equivariant morphism $\phi \colon X \to Z$, locally trivial in the Zariski topology, such that:*
  (i) *the fiber $F$ is a product $J_1 \times \cdots \times J_r$;*
 (ii) *each $J_i$ is the projective join of invariant subvarieties, and has Picard number 1;*
(iii) *each line in $X$ which intersects the open orbit of the torus, is contained in some fiber of $\phi$.*

Recall that a toric variety is Q-factorial when it has at worst quotient singularities. In example 5.7 we describe explicitly all varieties as in Theorem 4.2 up to dimension three.

An important ingredient in the proof is the description of $X$ when $\operatorname{Pic} X$ has rank one. In Theorem 2.6 we show that, in this case, $X$ is a projective join of invariant subvarieties. The existence of the fibration is then derived by a result in [BCD07].

Theorem 4.2 has an immediate application in the study of defective toric varieties.

---

2000 *Mathematics Subject Classification.* 14M25, 52B20.
This research has been partially supported by a Göran Gustafsson Stistelse Grant (2004-2005) and by the italian research project "Geometria sulle varietà algebriche" (COFIN 2004).

The dual defect of an embedded variety $X \subset \mathbb{P}^N$ is $\operatorname{def}(X) = N - 1 - \dim X^*$, where $X^* \subset (\mathbb{P}^N)^*$ is the dual variety of $X$. The embedding is said to be defective when $\operatorname{def}(X) > 0$, namely when $X^*$ is not a hypersurface. Varieties with positive dual defect are covered by lines and hence uniruled. In the smooth case they are quite exceptional, and much work has been devoted to their classification. In the singular case very little is known.

We use our results to compute the dual defect of $X \subset \mathbb{P}^N$ toric and $\mathbb{Q}$-factorial, generalizing the classification in [DR06] for the smooth case.

If $X$ is covered by lines and has Picard number 1, Theorem 2.6 implies that $X$ is a projective join and hence its dual defect is always positive. In fact $\operatorname{def}(X) = r - 1$, where $r$ is the maximal number of pairwise disjoint invariant subvarieties whose join gives the variety $X$ (see Corollary 2.8).

In the general case, $X$ inherits the structure described in Theorem 4.2. The dual defect of $X$ is related to the dual defect of the fiber $F$. In Theorem 5.2 we show that:
$$\operatorname{def}(X) = \max(0, \operatorname{def}(F) - \dim Z).$$

This yields an explicit formula for $\operatorname{def}(X)$:

**Theorem 5.2 - Corollary 5.4.** *Notation as above.*

$\operatorname{def}(F) = \max(0, \dim J_1 + \operatorname{def}(J_1) - \dim F, \ldots, \dim J_r + \operatorname{def}(J_r) - \dim F)$,

$\operatorname{def}(X) = \max(0, \dim J_1 + \operatorname{def}(J_1) - \dim X, \ldots, \dim J_r + \operatorname{def}(J_r) - \dim X)$.

*Moreover if $\operatorname{def}(X) > 0$, then there exists a morphism $\phi^* \colon X^*_{reg} \longrightarrow Z$ such that for a general point $z \in Z$, the closure of $(\phi^*)^{-1}(z)$ in $X^*$ is $(\phi^{-1}(z))^*$.*

Using the above results we are able to characterize defective $\mathbb{Q}$-factorial toric varieties as follows.

**Corollary 5.5.** *Let $X \subset \mathbb{P}^N$ be a toric and $\mathbb{Q}$-factorial variety. Then $X$ has positive dual defect if and only if there exists an elementary extremal contraction of fiber type $\psi \colon X \to Y$ whose general fiber is a projective join $J$ with $\dim J + \operatorname{def}(J) > \dim X$.*

As customary in toric geometry, this problem has a geometrical as well as a combinatorial significance. The theory of projective duality has attracted considerable attention in combinatorics, because it naturally arises in the study of $A$-discriminants, where $A$ is a finite subset of $\mathbb{Z}^n$. We refer to [GKZ94, DS02] for details. Let $X_A \subset \mathbb{P}^{|A|-1}$ be the variety associated to $A$. Then the $A$-discriminant $D_A$ is equal to a defining homogeneous polynomial of the dual variety $X_A^*$ if $\operatorname{def}(X_A) = 0$, and it is set to be 1 otherwise. The classification in the case $D_A = 1$ is an interesting and open problem. In [DS02, CC07] the case of codimension at most four is described.

Our results give an answer to this problem in any codimension, under suitable assumptions on the polytope $\operatorname{Conv}(A)$. The geometric description of a fibration arising from $\mathbb{Q}$-factorial toric varieties with dual defect yields a characterization of $\mathrm{Conv}(A)$ when $D_A = 1$. We show that when $D_A = 1$, $\mathrm{Conv}(A)$ is a $\pi$-*twisted Cayley sum* (see Definition 3.5), as follows.

**Theorem 6.1.** *Suppose that* $\mathrm{Conv}(A)$ *is a simple polytope of dimension* $n$.

*For every vertex* $v$ *of* $\mathrm{Conv}(A)$, *let* $P_v$ *be the translated polytope* $\mathrm{Conv}(A) - v$, *and let* $\eta_v$ *be the cone over* $P_v$ *in* $\mathbb{Q}^n$. *Assume that the semigroup* $\eta_v \cap \mathbb{Z}^n$ *is generated by* $A - v$.

*Then* $D_A = 1$ *if and only if there exists a surjection of lattices* $\pi \colon \mathbb{Z}^n \to \Lambda$ *such that:*

(i) $S := \mathrm{Conv}(\pi(A)) \subset \Lambda_{\mathbb{Q}}$ *is a lattice simplex* $\mathrm{Conv}(v_0, \ldots, v_{\dim S})$ *with* $\dim S + \mathrm{def}(S) > n$;

(ii) $A \subset \mathrm{Conv}(A_0, \ldots, A_{\dim S})$, *where* $A_i := A \cap \pi^{-1}(v_i)$;

(iii) *the polytopes* $\mathrm{Conv}(A_i) \subset \mathbb{Q}^n$ *are strictly combinatorially isomorphic for* $i = 0, \ldots, \dim S$.

See Definitions 2.5 and 1.4 for the notions of defect of a lattice simplex and of strictly combinatorially isomorphic lattice polytopes.

## Contents

## 1. Notation and preliminaries

A variety is an integral scheme, and we always work over the complex numbers. We say that a variety $X$ is $\mathbb{Q}$-factorial if every Weil divisor in $X$ has a multiple which is Cartier. The *Picard number* $\rho_X$ of $X$ is the rank of the group of Cartier divisors in $X$ modulo numerical equivalence. Throughout the paper, $n$ always denotes the dimension of $X$. We denote by $X_{reg}$ the smooth locus of $X$.

**The dual variety.** Let $X \subset \mathbb{P}^N$ be a variety. For any $p \in X_{reg}$, we denote by $\mathbb{T}_p X \subset \mathbb{P}^N$ the projective tangent space to $X$ at $p$.

If $q \in (\mathbb{P}^N)^*$, let $H_q \subset \mathbb{P}^N$ be the corresponding hyperplane. Set $\mathcal{J}_X^0 := \{(p,q) \in X_{reg} \times (\mathbb{P}^N)^* \mid H_q \supseteq \mathbb{T}_p X\}$ and let $\mathcal{J}_X$ be the closure of $\mathcal{J}_X^0$ in $X \times (\mathbb{P}^N)^*$. Then the dual variety $X^* \subset (\mathbb{P}^N)^*$ is the image of $\mathcal{J}_X$ under

the second projection:

(1)
$$\begin{array}{c} \mathcal{J}_X \subset \mathbb{P}^N \times (\mathbb{P}^N)^* \\ {}^{\pi_1}\swarrow \qquad \searrow^{\pi_2} \\ X \subset \mathbb{P}^N \qquad\qquad X^* \subset (\mathbb{P}^N)^* \end{array}$$

Observe that one can dually define $\mathcal{J}^0_{X^*}$ and $\mathcal{J}_{X^*}$, and that $\mathcal{J}_{X^*} = \mathcal{J}_X$ by reflexivity.

For any $q \in X^*$, we call $\pi_1(\pi_2^{-1}(q))$ the *contact locus* of $H_q$ with $X$. If $q$ is a general point, then the contact locus is a linear subspace of dimension $\mathrm{def}(X)$. In particular, if $X$ has positive dual defect, then it is covered by lines.

Notice that the dual defect can vary if we consider different embeddings of the same abstract variety, but it only depends on $X$ and $\mathcal{O}_X(1)$.

*Remark* 1.1. Let $X \subset \mathbb{P}^N$ and let $X' \subset \mathbb{P}^{N'}$ be the image of $X$ under the embedding defined by the complete linear system $|\mathcal{O}_X(1)|$. Then $\mathrm{def}(X) = \mathrm{def}(X')$ (see for instance [DR06, Proposition 1]).

**Projective joins.** Consider $Y_0, \ldots, Y_r \subset \mathbb{P}^N$ disjoint varieties such that for every $y_0 \in Y_0, \ldots, y_r \in Y_r$ the linear span $\overline{y_0 \cdots y_r}$ of the points $y_0, \ldots, y_r$ has dimension $r$. The *projective join* of $Y_0, \ldots, Y_r$ is defined as

$$J = J(Y_0, \ldots, Y_r) = \bigcup_{y_0 \in Y_0, \ldots, y_r \in Y_r} \overline{y_0 \cdots y_r}.$$

It is an irreducible variety of dimension equal to $\dim Y_0 + \cdots + \dim Y_r + r$. By definition $J$ is covered by lines. Moreover, it is well known that $J$ has dual defect $\mathrm{def}(J) \geq r$ (see the proof of Corollary 2.8).

**Notation for lattices and convex geometry.** For any lattice $\Lambda$, we set $\Lambda_{\mathbb{Q}} := \Lambda \otimes_{\mathbb{Z}} \mathbb{Q}$ and $\Lambda^{\vee} := \mathrm{Hom}_{\mathbb{Z}}(\Lambda, \mathbb{Z})$. If $\varphi \colon \Lambda \to \Gamma$ is a homomorphism of lattices, we denote by $\varphi_{\mathbb{Q}} \colon \Lambda_{\mathbb{Q}} \to \Gamma_{\mathbb{Q}}$ the induced $\mathbb{Q}$-homomorphism.

If $v_1, \ldots, v_m \in \Lambda_{\mathbb{Q}}$, we denote by $\langle v_1, \ldots, v_m \rangle$ the convex cone in $\Lambda_{\mathbb{Q}}$ generated by $v_1, \ldots, v_m$. For any subset $A \subset \Lambda_{\mathbb{Q}}$:

- $\mathrm{Conv}(A)$ is the minimal convex set containing $A$;
- $\mathrm{Aff}(A)$ is the minimal affine subspace of $\Lambda_{\mathbb{Q}}$ containing $A$;
- $\mathrm{RelInt}(A)$ is the interior of $A$ in $\mathrm{Aff}(A)$;
- $A^{\perp} = \{\varphi \in \Lambda_{\mathbb{Q}}^{\vee} \mid \varphi \text{ is trivial on } A\}$.

A *lattice polytope* in $\Lambda_{\mathbb{Q}}$ is a polytope whose vertices belong to $\Lambda$.

A *simplex* is a polytope of dimension $n$ with $n+1$ vertices.

A polytope of dimension $n$ is *simple* if every vertex is contained in exactly $n$ edges.

If $S$ is a finite set, we denote by $|S|$ its cardinality.

**Preliminaries in toric geometry.** We recall some notions and properties which will be useful in the sequel. We refer the reader to [Ful93, Ewa96] for more details.

A toric variety is a normal algebraic variety endowed with an (algebraic) action of a torus $T = (\mathbb{C}^*)^n$, and containing a dense open orbit.

Let $X$ be a toric variety of dimension $n$. The fan $\Sigma_X$ of $X$ is contained in the $\mathbb{Q}$-vector space $N_\mathbb{Q}$, where $N = \mathrm{Hom}(\mathbb{C}^*, T)$ is naturally identified with the group of one parameter subgroups of $T$. Moreover, $M$ is the dual lattice and $M_\mathbb{Q}$ the dual vector space. We denote by $(\ ,\ )$ the standard pairing on $N_\mathbb{Q} \times M_\mathbb{Q}$.

We say that a subset of $X$ is *invariant* if it is closed with respect to the action of $T$.

For any $\tau \in \Sigma_X$, denote by $O_\tau$ the orbit corresponding to $\tau$. The variety $V(\tau) := \overline{O_\tau}$ is an irreducible invariant subvariety of $X$. This gives a bijection among $(n-k)$-dimensional cones in $\Sigma_X$ and $k$-dimensional irreducible invariant subvarieties of $X$.

Moreover, we denote by $U_\tau$ the invariant affine open subset of $X$ corresponding to $\tau$. In particular, $U_0$ is the open orbit.

Let $G_X = \{x \in N \mid \langle x \rangle \in \Sigma_X \text{ and } x \text{ is primitive in } \langle x \rangle \cap N\}$, so that $G_X$ is in bijection with one dimensional cones in $\Sigma_X$. The invariant prime divisors are $D_x := V(\langle x \rangle)$ for $x \in G_X$.

When $X$ is $\mathbb{Q}$-factorial, we have $\rho_X = |G_X| - n$.

For any $u \in M$, we denote by $\chi^u$ the associated rational function on $X$.

**Toric varieties and polytopes.** Let $X$ be a toric variety with fan $\Sigma_X$, and let $L \in \mathrm{Pic}\, X$ be an ample line bundle. We denote by $\mathcal{P}_{(X,L)} \subset M_\mathbb{Q}$ the lattice polytope associated to the pair $(X, L)$. When $X \subset \mathbb{P}^N$, we just write $\mathcal{P}_X$ for $\mathcal{P}_{(X, \mathcal{O}_X(1))}$.

The polytope $\mathcal{P}_{(X,L)}$ is determined up to translation by elements of $M$. If $D = \sum_{x \in G_X} a_x D_x$ is an invariant Cartier divisor such that $L = \mathcal{O}_X(D)$, then the polytope is defined explicitly as

$$(2) \qquad \mathcal{P}_{(X,L)} = \{u \in M_\mathbb{Q} \mid (u, x) \geq -a_x \text{ for every } x \in G_X\}.$$

The following facts will be useful:
- there is a bijection among $k$-dimensional irreducible invariant subvarieties of $X$ and $k$-dimensional faces of $\mathcal{P}_{(X,L)}$;
- $X$ is $\mathbb{Q}$-factorial if and only if $\mathcal{P}_{(X,L)}$ is a simple polytope;
- $H^0(X, L) = \bigoplus_{u \in \mathcal{P}_{(X,L)} \cap M} \mathbb{C}\chi^u$.

*Remark* 1.2 ([Ful93], §1.5). Let $Z = V(\sigma)$ be an invariant subvariety of $X$, and $Q$ be the corresponding face of $\mathcal{P}_{(X,L)}$. Then:
- there exists $u \in M$ such that $\mathrm{Aff}(Q) + u = \sigma^\perp$;
- $Q + u \subset \sigma^\perp$ is a polytope associated to the pair $(Z, L_{|Z})$ (with respect to the lattice $M \cap \sigma^\perp$).

The very ampleness of $L$ translates to properties of the polytope $\mathcal{P}_{(X,L)}$.

**Lemma 1.3** ([Ful93], Lemma on p. 69). *For any vertex $v$ of $\mathcal{P}_{(X,L)}$, let $P_v$ be the translated polytope $\mathcal{P}_{(X,L)} - v$, and let $\eta_v$ be the cone over $P_v$.*

*Then $L$ is very ample if and only if for every vertex $v$ of $\mathcal{P}_{(X,L)}$ the semigroup $\eta_v \cap M$ is generated by $P_v \cap M$.*

**Definition 1.4** ([Ewa96], IV.2.11). *Two lattice polytopes $P$ and $Q$ in $M_{\mathbb{Q}}$ are* strictly combinatorially isomorphic *if there exists a bijection $\varphi$ from the set of faces of $P$ to the set of faces of $Q$, such that:*

(i) *$\varphi$ preserves inclusions*
(ii) *for every face $K$ of $P$ the affine subspaces $\mathrm{Aff}(K)$ and $\mathrm{Aff}(\varphi(K))$ are translates in $M_{\mathbb{Q}}$.*

By [Ewa96, IV.2.12], $P$ and $Q$ are strictly combinatorially isomorphic if and only if they define the same fan in $N_{\mathbb{Q}}$, so that $P = \mathcal{P}_{(X,L)}$ and $Q = \mathcal{P}_{(X,M)}$ with $L, M$ ample line bundles on $X$.

## 2. The case of Picard number one

Throughout this section, $X$ is a $\mathbb{Q}$-factorial projective toric variety with $\dim X = n$ and $\rho_X = 1$. We denote by $e_0, \ldots, e_n$ the elements of $G_X$. For every $i = 0, \ldots, n$ let $a_i \in \mathbb{Z}_{>0}$ be the least positive integer such that $a_i D_{e_i}$ is Cartier. We denote by $H$ the ample generator of $\mathrm{Pic}\, X \cong \mathbb{Z}$.

The following illustrates two elementary properties which will be needed in the sequel.

**Lemma 2.1.** *Let $C_0$ be a curve whose numerical class generates $\mathcal{N}_1(X) \cong \mathbb{Z}$. Then:*

(i) *$H \cong \mathcal{O}_X(a_i D_{e_i})$ for every $i = 0, \ldots, n$;*
(ii) *$H \cdot C_0 = 1$;*
(iii) *$\sum_{i=0}^n \frac{1}{a_i} e_i = 0$.*

*Proof.* Let $r := H \cdot C_0$ and $i \in \{0, \ldots, n\}$. Because $\mathcal{O}_X(a_i D_{e_i}) \in \mathrm{Pic}\, X$, there exists a unique $d_i \in \mathbb{Z}_{>0}$ such that $\mathcal{O}_X(a_i D_{e_i}) \cong H^{\otimes d_i}$. Moreover the minimality of $a_i$ gives $(a_i, d_i) = 1$.

Let $C \subset X$ be an invariant curve which is not contained in $D_{e_i}$. We have that $C \equiv mC_0$ for some $m \in \mathbb{Z}_{>0}$ and
$$C \cdot D_{e_i} = mC_0 \cdot D_{e_i} = \frac{m}{a_i} C_0 \cdot H^{\otimes d_i} = \frac{md_i r}{a_i}.$$

But, since $C \not\subseteq D_{e_i}$, we know that $\frac{1}{C \cdot D_{e_i}}$ must be an integer (see [Wiś02, §2]). It follows that $\frac{a_i}{md_i r} \in \mathbb{Z}$, $(d_i, a_i) = 1$, and hence $d_i = 1$. This proves (i).

Observe that by Reid's Toric Cone Theorem (see [Wiś02, §1]), $\mathcal{N}_1(X)$ is generated by classes of invariant curves. Therefore there is at least one invariant curve $C'$ which is numerically equivalent to $C_0$. Let $D_{e_{i_0}}$ be an invariant divisor which does not contain $C'$. For $b := \frac{a_{i_0}}{r}$ we have $C' \cdot D_{e_{i_0}} =$ $\frac{1}{b}$, thus $b \in \mathbb{Z}$ and $bD_{e_{i_0}}$ is Cartier. By the minimality of $a_{i_0}$, it must be $b = a_{i_0}$ and $r = 1$. This proves part $(ii)$.

Finally, the equality $C_0 \cdot D_{e_i} = \frac{1}{a_i}$ for any $i = 0, \ldots, n$ yields $(iii)$, see [Wiś02, §1] for details. ∎

**Lemma 2.2.** *Let $v \in N$ and let $\sigma \in \Sigma_X$ be such that $v \in \mathrm{RelInt}(\sigma)$. Let $I \subset \{0, \ldots, n\}$ be such that $\sigma = \langle e_i \rangle_{i \in I}$. Then*

$$v = \sum_{i \in I} \frac{m_i}{a_i} e_i, \quad \text{with } m_i \in \mathbb{Z}_{>0}.$$

*Proof.* Up to renumbering $e_0, \ldots, e_n$, we can assume that $I = \{0, \ldots, r\}$ with $r \in \{0, \ldots, n-1\}$. We know that $v = \sum_{i=0}^{r} b_i e_i$, with $b_i \in \mathbb{Q}_{>0}$. We have to show that $a_i b_i$ is an integer for every $i = 0, \ldots, r$.

Let $\eta := \langle e_0, \ldots, e_{n-1} \rangle$ and consider the affine invariant open subset $U_\eta$.

For a fixed $i \in \{0, \ldots, r\}$, $(a_i D_{e_i})_{|U_\eta}$ is Cartier and invariant. This implies that $(a_i D_{e_i})_{|U_\eta}$ is principal, namely there exists $\varphi_i \in M$ such that $\chi^{\varphi_i}$ is regular on $U_\eta$ and $(a_i D_{e_i})_{|U_\eta} = \mathrm{div}(\chi^{\varphi_i})$. It follows that $(e_i, \varphi_i) = a_i$ and $(e_j, \varphi_i) = 0$ for $j = 0, \ldots, n-1$, $j \neq i$ (see [Ful93, §3.3]). Then we have $(v, \varphi_i) = (b_i e_i, \varphi_i) = a_i b_i \in \mathbb{Z}$. ∎

**2.1. Degrees of orbits of one parameter subgroups.** Recall that $N$ is canonically identified with the group of one parameter subgroups of the torus $T$. For any $v \in N$, we denote by $\lambda_v \colon \mathbb{C}^* \to T$ the associated one parameter subgroup.

Let $v$ be a non zero element in $N$. Since $X$ is complete, there are unique cones $\sigma, \sigma' \in \Sigma_X$ such that $v \in \mathrm{RelInt}(\sigma)$ and $-v \in \mathrm{RelInt}(\sigma')$. It follows that for any $x \in U_0$ we have (see [Ful93, §2.3])

$$(3) \qquad \lim_{z \to 0} \lambda_v(z) \cdot x \in O_\sigma \quad \text{and} \quad \lim_{z \to \infty} \lambda_v(z) \cdot x \in O_{\sigma'}.$$

Let $x \in U_0$ and let $C_v$ be the closure in $X$ of the orbit $\lambda_v(\mathbb{C}^*) \cdot x$. Then $C_v$ is a rational curve in $X$ which intersects $X \smallsetminus U_0$ exactly in two distinct points. If we vary the point $x \in U_0$, the corresponding curves $C_v$ are translated under the torus action, in particular they are algebraically equivalent.

Observe also that $C_v = C_u$ for any $u \in (\mathbb{Q} \cdot v) \cap N$, $u \neq 0$. If $v$ is primitive in $N$, then the map $\mathbb{P}^1 \to C_v$ induced by $\lambda_v$ is birational (and bijective).

**Lemma 2.3.** *Let $v \in N$ be a non zero primitive element. If $v = \sum_{i \in I} \frac{m_i}{a_i} e_i$ with $I \subsetneq \{0, \ldots, n\}$ and $m_i \in \mathbb{Z}_{>0}$ for all $i \in I$, then*

$$C_v \cdot H = \max_{i \in I}(m_i),$$

*where $H$ is the ample generator of $\mathrm{Pic}\, X$.*

*Proof.* Up to renumbering $e_0, \ldots, e_n$, we can assume that $I = \{0, \ldots, r\}$ with $r \in \{0, \ldots, n-1\}$, and that $m_0$ is the maximum among $m_0, \ldots, m_r$.

Let $\sigma, \sigma' \in \Sigma_X$ be such that $v \in \operatorname{RelInt}(\sigma)$ and $-v \in \operatorname{RelInt}(\sigma')$. Then $\sigma = \langle e_0, \ldots, e_r \rangle$, and using Lemma 2.1 $(iii)$ we have that

$$-v = m_0 \left( \sum_{i=0}^{n} \frac{1}{a_i} e_i \right) - \sum_{i=0}^{r} \frac{m_i}{a_i} e_i = \sum_{i=0}^{r} \frac{m_0 - m_i}{a_i} e_i + \sum_{i=r+1}^{n} \frac{m_0}{a_i} e_i.$$

Because the coefficients of $e_1, \ldots, e_n$ are non negative and the coefficient of $e_0$ is zero it follows that $e_0 \notin \sigma'$.

Let $x_0 \in U_0$ be such that $C_v = \overline{\lambda_v(\mathbb{C}^*) \cdot x_0}$, and let $p := \lim_{z \to 0} \lambda_v(z) \cdot x_0$ and $q := \lim_{z \to \infty} \lambda_v(z) \cdot x_0$. By (3) we know that $p \in D_{e_0}$ and $q \notin D_{e_0}$, so $C_v \cap D_{e_0} = \{p\}$. Observe that $H \cong \mathcal{O}_X(a_0 D_{e_0})$ by Lemma 2.1 $(i)$.

Since $v$ is primitive in $N$, $\lambda_v$ induces a birational, bijective map $\psi \colon \mathbb{P}^1 \to C_v$ such that $\psi(0) = p$ and $\psi(\infty) = q$.

Let $\eta := \langle e_0, \ldots, e_{n-1} \rangle$ and consider the affine open subset $U_\eta$. Observe that $C_0 := C_v \cap U_\eta$ is exactly $C_v \smallsetminus \{q\}$ and that $\psi(\mathbb{A}_z^1) = C_0$.

Since $(a_0 D_{e_0})_{|U_\eta}$ is Cartier and invariant, there exists $\varphi \in M$ such that $\chi^\varphi$ is regular on $U_\eta$ and $(a_0 D_{e_0})_{|U_\eta} = \operatorname{div}(\chi^\varphi)$. It follows that $(e_0, \varphi) = a_0$, $(e_i, \varphi) = 0$ for $i = 1, \ldots, n-1$ (see [Ful93, §3.3]) and thus

$$(v, \varphi) = \left( \frac{m_0}{a_0} e_0, \varphi \right) = m_0.$$

The pull-back under $\psi$ of the regular function $\chi^\varphi$ to $\mathbb{A}_z^1$ is $z^{(v,\varphi)} = z^{m_0}$ (see [Ful93, §2.3]), and therefore the degree of the line bundle $\psi^*(H)$ on $\mathbb{P}^1$ is $m_0$. Since $\psi$ is birational, we conclude that $C_v \cdot H = m_0$. ∎

**Corollary 2.4.** *The following statements are equivalent:*
  *(i)* *there exists a non zero $v \in N$ such that $C_v \cdot H = 1$;*
  *(ii)* *there exists a proper subset $I \subsetneq \{0, \ldots, n\}$ such that $\sum_{i \in I} \frac{1}{a_i} e_i \in N$.*

*Proof.* Assume $(i)$. Let $u$ be a generator of the group $(\mathbb{Q} \cdot v) \cap N$. Then $C_u = C_v$, and Lemma 2.3 implies that there exists a subset $I \subset \{0, \ldots, n\}$ such that $u = \sum_{i \in I} \frac{1}{a_i} e_i$.

Conversely, assume $(ii)$ and let $v := \sum_{i \in I} \frac{1}{a_i} e_i$. Lemma 2.2 implies that $v$ is primitive and Lemma 2.3 yields $C_v \cdot H = 1$. ∎

**2.2. Projective joins of toric varieties.** We are now going to use the results of 2.1 to show that every $\mathbb{Q}$-factorial toric variety $X \subset \mathbb{P}^N$ with Picard number one and covered by lines is a projective join. The crucial step is to find a line in $X$ which intersects $U_0$, and is invariant with respect to some one-parameter subgroup of $T$.

We will also relate the dual defect of $X$ to an invariant of the lattice simplex $\mathcal{P}_X$, that we call the lattice defect of $\mathcal{P}_X$.

**Definition 2.5.** Let $\Lambda$ be a lattice and $P \subset \Lambda_\mathbb{Q}$ a lattice simplex.

The *lattice defect* $\operatorname{def}(P)$ of $P$ is the maximal $d \in \mathbb{Z}_{\geq 0}$ such that there exist $d + 1$ faces $Q_0, \ldots, Q_d$ of $P$ with the following properties:
  $(i)$ $Q_0, \ldots, Q_d$ are pairwise disjoint;

(ii) $P \cap \Lambda = \bigcup_{i=0}^{d}(Q_i \cap \Lambda)$.

Notice that $\mathrm{def}(P) \leq \dim P$, and equality holds if and only if $P$ is isomorphic (as a lattice polytope) to the standard simplex. See example 6.3 for examples of 2-dimensional simplices with lattice defect 0, 1, and 2.

**Theorem 2.6.** *Let $X \subset \mathbb{P}^N$ be a $\mathbb{Q}$-factorial toric variety with $\rho_X = 1$. The following are equivalent:*

(i) *$X$ is covered by lines;*
(ii) *$\mathcal{O}_X(1)$ is a generator of $\mathrm{Pic}\, X$, and there exists a proper subset $I \subsetneq \{0, \ldots, n\}$ such that $\sum_{i \in I} \frac{1}{a_i} e_i \in N$;*
(iii) *there exists a line $L \subseteq X$ such that $|L \cap (X \smallsetminus U_0)| = 2$;*
(iv) *$X$ is the projective join of two disjoint invariant subvarieties;*
(v) *$\mathrm{def}(X) > 0$;*
(vi) *$\mathrm{def}(\mathcal{P}_X) > 0$.*

*Proof.* $(i) \Rightarrow (ii)$ Since $X$ contains a line, Lemma 2.1 implies that $X$ is embedded by the ample generator $H = \mathcal{O}_X(1)$ of $\mathrm{Pic}\, X$.

Let $L \subset X$ be a line such that $L \cap U_0 \neq \emptyset$. Clearly $L \cong \mathbb{P}^1$ can not be contained in $U_0$ which is affine. Let $p \in L \cap (X \smallsetminus U_0)$ and $\sigma = \langle e_i \rangle_{i \in I} \in \Sigma_X$ be such that $p \in O_\sigma$.

Consider an integer $m \in \mathbb{Z}_{>0}$ such that

$$v := m \sum_{i \in I} \frac{1}{a_i} e_i \in N.$$

If for a point $x \in U_0$ the curve $\overline{\lambda_v(\mathbb{C}^*) \cdot x}$ is a line, then Lemma 2.3 implies that $\sum_{i \in I} \frac{1}{a_i} e_i \in N$ and thus $(ii)$.

Consider the morphism $\varphi \colon \mathbb{C}^* \times L \to X$ defined by $\varphi(z, x) = \lambda_v(z) \cdot x$. This gives a rational map $\varphi \colon \mathbb{P}^1 \times L \dashrightarrow X$. Let $S \to \mathbb{P}^1 \times L$ be a sequence of blow-ups giving a minimal resolution $\overline{\varphi} \colon S \to X$ of $\varphi$. This construction is described in detail in [Wiś02, p. 253].

$$\begin{array}{ccc}
 & S & \\
{}^{f}\swarrow & \downarrow & \searrow{}^{\overline{\varphi}} \\
\mathbb{P}^1 \longleftarrow & \mathbb{P}^1 \times L \dashrightarrow_{\varphi} & X
\end{array}$$

Let $f \colon S \to \mathbb{P}^1$ be the natural projection, and let $L' := \overline{\varphi}_*(f^{-1}(\infty))$. Since $L = \overline{\varphi}_*(f^{-1}(1))$, we see that the two curves are algebraically equivalent and therefore $L'$ is a line. Observe also that $L'$ is closed with respect to $\lambda_v$. It remains to show that $L' \cap U_0 \neq \emptyset$. This would imply $L' = \overline{\lambda_v(\mathbb{C}^*) \cdot x}$ for some $x \in U_0$, and hence conclude the proof.

By (3), $\lambda_v$ fixes each point in $O_\sigma$, in particular the point $p$, so $\varphi(\mathbb{C}^* \times p) = p$ and $p \in O_\sigma \cap L'$.

Now let $x_0 \in L \cap U_0$ and $\sigma' := \langle e_i \rangle_{i \notin I}$. We have that

$$-v = m \sum_{i \notin I} \frac{1}{a_i} e_i,$$

which by (3) implies

$$q := \lim_{z \to \infty} \lambda_v(z) \cdot x_0 \in O_{\sigma'} \cap L'.$$

Finally, notice that $p \in D_{e_i}$ if and only if $i \in I$, while $q \in D_{e_i}$ if and only if $i \notin I$. Hence there is no $D_{e_i}$ containing $L'$, and $L' \cap U_0 \neq \emptyset$.

$(ii) \Rightarrow (iii)$ By Corollary 2.4.

$(iii) \Rightarrow (iv)$ Let $L \cap (X \smallsetminus U_0) =: \{p, q\}$.

Since $X$ contains a line, Lemma 2.1 yields $\mathcal{O}_X(1) = H \cong \mathcal{O}_X(a_i D_{e_i})$ for any $i = 0, \ldots, n$. In particular $L$ must intersect each $D_{e_i}$ exactly in one point, which is either $p$ or $q$. So if we define $I_p := \{i \,|\, p \in D_{e_i}\}$ and $I_q := \{i \,|\, q \in D_{e_i}\}$, we have that $I_p \cap I_q = \emptyset$ and $I_p \cup I_q = \{0, \ldots, n\}$.

Call $V_p$ and $V_q$ the closures of the orbits of $p$ and $q$ respectively. Then

$$V_p = \bigcap_{i \in I_p} D_{e_i} \quad \text{and} \quad V_q = \bigcap_{i \in I_q} D_{e_i}, \quad \text{so} \quad V_p \cap V_q = \bigcap_{i=0}^{n} D_{e_i} = \emptyset,$$

and $\dim V_p + \dim V_q = (n - |I_p|) + (n - |I_q|) = n - 1$. Hence the projective join $J = J(V_p, V_q)$ of $V_p$ and $V_q$ is an irreducible variety of dimension

$$\dim J = \dim V_p + \dim V_q + 1 = n.$$

Since $X$ and $J$ are irreducible of the same dimension, $X \subseteq J$ would imply $X = J$.

Let $x, x' \in U_0$ be such that $x \in L$ and $x' = t \cdot x$ for some $t \in T$. Consider the curve $L' = t(L)$. It is algebraically equivalent to $L$, so it is a line, and clearly $x' \in L'$. Moreover $t \cdot p \in L' \cap V_p$, $t \cdot q \in L' \cap V_q$ and hence $L' \subset J$. In particular we have that $x' \in J$. This shows that the open subset $U_0$ is contained in $J$ and hence $X \subseteq J$.

$(iv) \Rightarrow (v)$ and $(v) \Rightarrow (i)$ : clear.

$(iv) \Rightarrow (vi)$ We can assume that $X \subset \mathbb{P}^N$ is not contained in a hyperplane and is embedded by the complete linear system $|\mathcal{O}_X(1)|$, so $N = h^0(\mathcal{O}_X(1)) - 1 = |\mathcal{P}_X \cap M| - 1$. We have $X = J(V_0, V_1)$, $V_0$ and $V_1$ disjoint invariant subvarieties. Consider the face $Q_i$ of $\mathcal{P}_X$ corresponding to $V_i$. Notice that $Q_0 \cap Q_1 = \emptyset$.

We have $V_i \subset \mathbb{P}^{N_i}$, where $N_i := |Q_i \cap M| - 1$ and $\mathbb{P}^{N_0} \cap \mathbb{P}^{N_1} = \emptyset$ in $\mathbb{P}^N$. Since $X = J(V_0, V_1)$ is not contained in a hyperplane, we get $N = N_0 + N_1 + 1$, namely $\mathcal{P}_X \cap M = (Q_0 \cap M) \cup (Q_1 \cap M)$. This implies that $\mathrm{def}(\mathcal{P}_X) > 0$.

$(vi) \Rightarrow (ii)$ Since $\mathrm{def}(\mathcal{P}_X) > 0$, there are two disjoint faces $Q_0, Q_1$ of $\mathcal{P}_X$ such that $\mathcal{P}_X \cap M = (Q_0 \cap M) \cup (Q_1 \cap M)$.

Let $r \geq 1$ be such that $\mathcal{O}_X(1) \cong H^{\otimes r}$. This yields $\mathcal{P}_X = r\mathcal{P}_{(X,H)}$. If $r \geq 2$, then every edge of $\mathcal{P}_X$ contains some lattice point in its relative interior. This gives a contradiction. Therefore $r = 1$ and $\mathcal{O}_X(1)$ is a generator of $\operatorname{Pic} X$.

Let $\varepsilon_1, \ldots, \varepsilon_n$ be the basis of $M_\mathbb{Q}$ dual to the basis $e_1, \ldots, e_n$ of $N_\mathbb{Q}$. Considering $a_0 D_0 \in |\mathcal{O}_X(1)|$ and using the explicit description (2) of $\mathcal{P}_X$, we see that $\mathcal{P}_X = \operatorname{Conv}(0, a_1 \varepsilon_1, \ldots, a_n \varepsilon_n)$.

Up to renumbering, we can assume that $Q_0 = \operatorname{Conv}(a_1 \varepsilon_1, \ldots, a_r \varepsilon_r)$ and $Q_1 = \operatorname{Conv}(0, a_{r+1} \varepsilon_{r+1}, \ldots, a_n \varepsilon_n)$ for some $r \in \{1, \ldots, n\}$.

We claim that $v_0 := \sum_{i=1}^r \frac{1}{a_i} e_i \in N$. Since $\mathcal{O}_X(1)$ is very ample, Lemma 1.3 implies that $\mathcal{P}_X$ contains a basis of $M$. It is then enough to prove that $(v_0, u) \in \mathbb{Z}$ for any $u \in \mathcal{P}_X \cap M$. Such $u$ must lie either in $Q_0$ or in $Q_1$. We have $(v_0, u) = 1$ in the first case and $(v_0, u) = 0$ in the second, so $v_0 \in N$. ∎

*Remark* 2.7. Let $X \subset \mathbb{P}^N$ be a toric variety and $V_0, \ldots, V_r$ pairwise disjoint invariant subvarieties. Then $\overline{x_0 \cdots x_r} = \mathbb{P}^r$ for every $x_0 \in V_0, \ldots, x_r \in V_r$.

It is enough to see this when $X$ is linearly normal. In this case, we have $V_i \subseteq \mathbb{P}^{N_i}$ with $\mathbb{P}^{N_0}, \ldots, \mathbb{P}^{N_r}$ pairwise disjoint in $\mathbb{P}^N$.

**Corollary 2.8.** *Let $X \subset \mathbb{P}^N$ be a ℚ-factorial toric variety with $\rho_X = 1$. Suppose that $X$ is covered by lines, and set $d := \operatorname{def}(X)$.*

*Then there are uniquely determined $d+1$ invariant subvarieties $V_0, \ldots, V_d \subset X$ such that:*
  (i) $V_0, \ldots, V_d$ are not covered by lines;
  (ii) $V_0, \ldots, V_d$ are pairwise disjoint;
  (iii) $X$ is the projective join of $V_0, \ldots, V_d$.

**Corollary 2.9.** *Let $X \subset \mathbb{P}^N$ be a ℚ-factorial toric variety with $\rho_X = 1$. Then $\operatorname{def}(X) = \operatorname{def}(\mathcal{P}_X)$.*

*Remark* 2.10. Recall that $n$-dimensional linear subspaces have dual defect $n$. Observe that, contrary to what happens in the smooth case, the defect can assume any value $d$ between $1$ and $n$. For instance, consider a (non linear) Veronese embedding $Y$ of $\mathbb{P}^{n-d}$. Then the cone over $Y$ with vertex a linear subspace of dimension $d-1$ has dual defect $d$.

In example 5.7 we list all cases for $n \leq 3$.

*Proof of Corollary 2.8.* Let's show, by induction on the dimension, that we can find pairwise disjoint invariant subvarieties, not covered by lines, such that $X$ is their projective join.

By Theorem 2.6, $X$ is the projective join of two disjoint invariant subvarieties $W_0$ and $W_1$. Observe that $W_0$ and $W_1$ are toric, ℚ-factorial, and have Picard number 1. By induction, $W_i$ is the projective join of $r_i \in \mathbb{Z}_{\geq 1}$ pairwise disjoint invariant subvarieties $W_{i1}, \ldots, W_{ir_i}$, which are not covered by lines. Now $W_{ij}$ is also an invariant subvariety of $X$, and $X$ is the projective join of the $W_{ij}$'s for $i = 0, 1$ and $j = 1, \ldots, r_i$.

For simplicity rename $W_{01}, \ldots, W_{0r_0}, W_{11}, \ldots, W_{1r_1}$ as $V_0, \ldots, V_r$. Let's show that $r = \operatorname{def}(X)$. Choose general points $p_i \in V_i$ and let $L := \overline{p_0 \cdots p_r}$ be their linear span. By Terracini's Lemma, the projective tangent space to $X$ is constant along $L$. It follows that if a hyperplane $H$ is tangent to $X$ at a general point of $L$, then $L$ is entirely contained in the contact locus of $H$ with $X$.

Choose such an $H$, and let $P = \mathbb{P}^{\operatorname{def} X}$ be its contact locus. We have that $P \supseteq L$.

Lemma 2.1 implies that $\mathcal{O}_X(1) = H \cong \mathcal{O}_X(a_j D_j)$ for any $j = 0, \ldots, n$. Hence each $D_j$ is (set theoretically) a hyperplane section of $X$. Each invariant subvariety of $X$ is (set theoretically) a complete intersection of the $D_j$'s. Then for each $i = 0, \ldots, r$ there exists a linear subspace $T_i \subset \mathbb{P}^N$ such that $V_i = X \cap T_i$. It follows that $P \cap V_i = P \cap T_i$ is a linear subspace contained in $V_i$. Because there are no lines in $V_i$ through a general point it must be $P \cap V_i = \{p_i\}$.

Suppose that $P \supsetneq L$, and let $q$ be a general point in $P \smallsetminus L$. Then $q$ lies on some $L' = \overline{p'_0 \cdots p'_r}$ which again must be contained in $P$. But then $p'_i = p_i$ for all $i$, which is a contradiction. We conclude that $P = L$ and $r = \operatorname{def}(X)$.

Finally, observe that any invariant subvariety $V' \subset X$ which is not covered by lines must be contained in some $V_i$. This yields uniqueness for $V_0, \ldots, V_{\operatorname{def}(X)}$. ∎

**Corollary 2.11.** *Let $X \subset \mathbb{P}^N$ be a $\mathbb{Q}$-factorial toric variety with $\rho_X = 1$. Suppose that $\mathcal{O}_X(1)$ is a generator of $\operatorname{Pic} X$, and that $D_{e_0}, \ldots, D_{e_r}$ are Cartier, with $r \in \{0, \ldots, n\}$. For any $i = 0, \ldots, r$ let $p_i$ be the unique fixed point not lying in $D_{e_i}$.*

*Then $X$ is the cone over $D_{e_0} \cap \cdots \cap D_{e_r}$ with vertex the linear span of $p_0, \ldots, p_r$.*

## 3. Toric fibrations

In this section we describe toric fibrations and collect some results which will be used in the sequel.

Throughout the section, $X$ is a $\mathbb{Q}$-factorial, projective toric variety of dimension $n$.

**Definition 3.1.** A *toric fibration* on $X$ is a surjective, flat, equivariant morphism $f\colon X \to Y$ with connected fibers, where $Y$ is a projective toric variety with $\dim Y < n$. A toric fibration is called *elementary* if $\rho_X - \rho_Y = 1$.

There is a simple combinatorial characterization of toric fibrations in terms of the fan of $X$.

**Lemma 3.2.** *Let $\Delta$ be a sublattice of $N$. The inclusion $\Delta \hookrightarrow N$ induces a toric fibration $f\colon X \to Y$ if and only if the following conditions are satisfied:*
  (i) $\Delta_{\mathbb{Q}} \cap N = \Delta$, namely $\Delta$ is a primitive sublattice of $N$;
  (ii) *for every $n$-dimensional $\sigma \in \Sigma_X$ we have $\sigma = \tau + \eta$ with $\tau, \eta \in \Sigma_X$, $\tau \subset \Delta_{\mathbb{Q}}$, and $\eta \cap \Delta_{\mathbb{Q}} = \{0\}$.*

*Proof.* See [Rei83, Theorem 2.4] for the case of an elementary fibration, and [Ewa96, §VI.6] for the general case. ∎

*Remark* 3.3. The relation between $f$ and $\Delta$ is the following: the lattice $N/\Delta$ is canonically identified with the lattice of one parameter subgroups of the torus $T_Y$, and $f \colon X \to Y$ is induced by the projection $N \to N/\Delta$. Let $\Sigma_F := \{\sigma \in \Sigma_X \,|\, \sigma \subset \Delta_\mathbb{Q}\}$. The fan $\Sigma_F$ is the fan (with respect to $\Delta$) of a projective toric variety $F$, which is isomorphic to every fiber of $f$, with the reduced scheme structure. *Both $Y$ and $F$ are $\mathbb{Q}$-factorial,* and $\dim F = \operatorname{rk} \Delta$. The varieties $V(\tau)$ for $\tau \in \Sigma_F$, $\dim \tau = \dim F$ are the invariant sections of $f$, while the invariant fibers are $V(\eta)$ for any $\eta \in \Sigma_X$ such that $\dim \eta = \dim Y$ and $\eta \cap \Delta_\mathbb{Q} = \{0\}$.

*Example* 3.4. Consider the Hirzebruch surface $\mathbb{F}_1 = \mathbb{P}(\mathcal{O}_{\mathbb{P}^1} \oplus \mathcal{O}_{\mathbb{P}^1}(1))$.

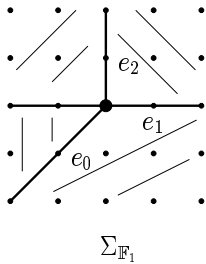

$\Sigma_{\mathbb{F}_1}$

The two-dimensional cones in the fan of $\mathbb{F}_1$ are $\langle e_1, e_2 \rangle$, $\langle -e_1, e_2 \rangle$, $\langle e_1, e_0 \rangle$, and $\langle -e_1, e_0 \rangle$, where $e_1, e_2$ is a basis of $N$ and $e_0 = -e_1 - e_2$. The fibration $\mathbb{F}_1 \to \mathbb{P}^1$ corresponds to the sublattice $\Delta = \mathbb{Z} \cdot e_1$.

When $X$ is polarized by an ample line bundle $L$, we call the fibration a *polarized toric fibration* and denote it by $(f \colon X \to Y, L)$. We give a dual combinatorial characterization of polarized toric fibrations in terms of the polytope $\mathcal{P}_{(X,L)}$. The following definition is a generalization of the Cayley sum of finitely many polytopes.

**Definition 3.5.** Let $\pi \colon M \to \Lambda$ be a surjective map of lattices and let $R_1, \ldots, R_l \subset M_\mathbb{Q}$ be lattice polytopes. Assume that:
- $\pi_\mathbb{Q}(R_i) = v_i \in \Lambda$ for every $i = 1, \ldots, l$;
- $v_1, \ldots, v_l$ are all distinct and are the vertices of the polytope $\operatorname{Conv}(v_1, \ldots, v_l) \subset \Lambda_\mathbb{Q}$;
- $R_1, \ldots, R_l$ are strictly combinatorially isomorphic.

Then we define the $\pi$-*twisted Cayley sum* of $R_1, \ldots, R_l$ to be the polytope
$$\mathcal{C}(R_1, \ldots, R_l, \pi) := \operatorname{Conv}(R_1, \ldots, R_l).$$

The proof of Lemma 3.6 will show that $\mathcal{C}(R_1, \ldots, R_l, \pi)$ has the same combinatorial type as the product $R_1 \times \operatorname{Conv}(v_1, \ldots, v_l)$. However in general they are not isomorphic as lattice polytopes, see example 3.7.

**Lemma 3.6.** *Let $L \in \operatorname{Pic} X$ be an ample line bundle, and $\Delta$ a sublattice of $N$. Set $\Lambda := \Delta^\vee$ and let $\pi \colon M \to \Lambda$ be the map dual to the inclusion* $j\colon \Delta \hookrightarrow N$. Then $j$ induces a polarized toric fibration $(f\colon X \to Y, L)$ if and only if:

(i) $\pi$ is surjective;

(ii) $\mathcal{P}_{(X,L)}$ is the $\pi$-twisted Cayley sum of some lattice polytopes $R_1, \ldots, R_l \subset M_{\mathbb{Q}}$.

When these conditions are fulfilled, we also have:
$$\pi_{\mathbb{Q}}(\mathcal{P}_{(X,L)}) = \mathcal{P}_{(F, L_{|F})}, \quad R_i - u_i = \mathcal{P}_{(Y_i, L_{|Y_i})}$$

where $F$ is a general fiber, $Y_1, \ldots, Y_l$ the invariant sections of $f$, and $u_i \in M$ is such that $\pi(u_i) = \pi(R_i)$.

*Example* 3.7. Consider the surface $\mathbb{F}_1$ as in example 3.4. Let $E \subset \mathbb{F}_1$ be the $(-1)$-curve, $F$ a fiber of the fibration onto $\mathbb{P}^1$, and $L = \mathcal{O}_{\mathbb{F}_1}(2F + E)$. Then $L$ is very ample and embeds $\mathbb{F}_1$ as a degree 3 surface $S \subset \mathbb{P}^4$. A polytope $P = \mathcal{P}_{(\mathbb{F}_1, L)}$ is described in the following figure, where $\varepsilon_1, \varepsilon_2$ is the basis of $M$ dual to $e_1, e_2$.

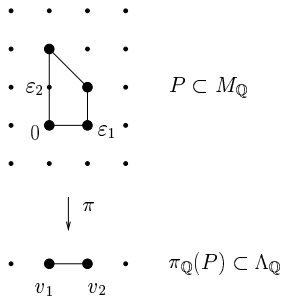

We have $R_i = P \cap \pi_{\mathbb{Q}}^{-1}(v_i)$ for $i = 1, 2$.

Recall that for any $a \in \mathbb{Z}_{\geq 1}$, the polytope $\mathcal{P}_{(\mathbb{P}^1, \mathcal{O}_{\mathbb{P}^1}(a))}$ is a one-dimensional lattice polytope containing $a + 1$ lattice points. Here, the fiber $F$ and the invariant section $E$ are embedded as lines in $\mathbb{P}^4$. In fact the corresponding polytopes, $\pi_{\mathbb{Q}}(P)$ and $R_2$, both contain 2 lattice points. The other invariant section is embedded as a conic and in fact the corresponding face, $R_1$, contains 3 lattice points.

*Proof.* Observe first of all that 3.2 $(i)$ is equivalent to 3.6 $(i)$.

Assume that $j \colon \Delta \hookrightarrow N$ induces a toric fibration $f \colon X \to Y$. Let $F$ be a general fiber, and $S := \mathcal{P}_{(F, L_{|F})} \subset \Lambda_{\mathbb{Q}}$. Denote by $v_1, \ldots, v_l$ the vertices of $S$. Every $v_i$ corresponds to a fixed point of $F$. Call $Y_i = V(\tau_i)$ the invariant section of $f$ passing through that point, with $\tau_i \in \Sigma_X$, $\dim \tau_i = \dim F$ and $\tau_i \subset \Delta_{\mathbb{Q}}$ (see remark 3.3). Finally, let $R_i$ be the face of $\mathcal{P}_{(X,L)}$ corresponding to $Y_i$.

Observe that $\mathrm{Aff}(\tau_i) = \Delta_{\mathbb{Q}}$, so that $\tau_i^{\perp} = \Delta_{\mathbb{Q}}^{\perp} = \ker \pi_{\mathbb{Q}}$. By remark 1.2, there exists $u_i \in M$ such that:

- $\mathrm{Aff}(R_i) + u_i = \ker \pi_{\mathbb{Q}}$;
- $R_i + u_i = \mathcal{P}_{(Y_i, L_{|Y_i})}$.

This says that $R_1, \ldots, R_l$ are strictly combinatorially isomorphic (because every $Y_i$ is isomorphic to $Y$), and that $\pi_\mathbb{Q}(R_i)$ is a point.

Since the $Y_i$'s are pairwise disjoint, the same holds for the $R_i$'s. If $s$ is the number of fixed points of $Y$, then each $R_i$ has $s$ vertices. On the other hand, we know that $F$ has $l$ fixed points, and $X$ must have $sl$ fixed points. So $\mathcal{P}_{(X,L)}$ has $sl$ vertices, namely the union of all vertices of the $R_i$'s. Therefore $\mathcal{P}_{(X,L)} = \mathrm{Conv}(R_1, \ldots, R_l) = \mathcal{C}(R_1, \ldots, R_l, \pi)$.

Let $D = \sum_{x \in G_X} a_x D_x$ be an invariant Cartier divisor on $X$ such that $L = \mathcal{O}_X(D)$. Since $F$ is a general fiber, we have $D_x \cap F \neq \emptyset$ if and only if $x \in \Delta$, and $D_{|F} = \sum_{x \in G_X \cap \Delta} a_x D_{x|F}$. Using the explicit description (2) for the polytopes $\mathcal{P}_{(X,L)}$ and $S$, it is easy to see that in fact $\pi_\mathbb{Q}(R_i) = v_i$ and $\pi_\mathbb{Q}(\mathcal{P}_{(X,L)}) = S$. This shows the last part of the statement.

Conversely, assume that $\mathcal{P}_{(X,L)} = \mathcal{C}(R_1, \ldots, R_l, \pi)$ where $R_1, \ldots, R_l$ are as in Definition 3.5.

Since $v_i$ is a vertex of $\pi_\mathbb{Q}(\mathcal{P}_{(X,L)})$, $R_i$ is a face of $\mathcal{P}_{(X,L)}$ for every $i = 1, \ldots, l$. Let $Y$ be the projective toric variety defined by the polytopes $R_i$. Observe that $\mathrm{Aff}(R_i)$ is a translate of $\ker \pi_\mathbb{Q}$, and $(\ker \pi)^\vee = N/\Delta$. So the fan $\Sigma_Y$ is contained in $(N/\Delta)_\mathbb{Q}$.

Let $\gamma \in \Sigma_Y$ and for every $i = 1, \ldots, l$ let $w_i$ be the vertex of $R_i$ corresponding to $\gamma$. Let's show that $Q := \mathrm{Conv}(w_1, \ldots, w_l)$ is a face of $\mathcal{P}_{(X,L)}$.

Observe first of all that $(\pi_\mathbb{Q})_{|\mathrm{Aff}(Q)} \colon \mathrm{Aff}(Q) \to \Lambda_\mathbb{Q}$ is bijective. Let $H$ be the linear subspace of $M_\mathbb{Q}$ which is a translate of $\mathrm{Aff}(Q)$. Then we have $M_\mathbb{Q} = H \oplus \ker \pi_\mathbb{Q}$. Dually $N_\mathbb{Q} = \Delta_\mathbb{Q} \oplus H^\perp$, where $H^\perp$ projects isomorphically onto $(N/\Delta)_\mathbb{Q}$.

Let $u \in H^\perp$ be such that its image in $(N/\Delta)_\mathbb{Q}$ is contained in the interior of $\gamma$. Then for every $i = 1, \ldots, l$ we have that (see [Ful93, §1.5]):

$$(u, x) \geq (u, w_i) \quad \text{for every } x \in R_i,$$
$$(u, x) = (u, w_i) \quad \text{if and only if } x = w_i.$$

Moreover $u$ is constant on $\mathrm{Aff}(Q)$, namely there exists $m_0 \in \mathbb{Q}$ such that $(u, z) = m_0$ for every $z \in Q$.

Any $z \in P$ can be written as $z = \sum_{i=1}^l \lambda_i z_i$, with $z_i \in R_i$, $\lambda_i \geq 0$ and $\sum_{i=1}^l \lambda_i = 1$. Then

$$(u, z) = \sum_{i=1}^l \lambda_i (u, z_i) \geq \sum_{i=1}^l \lambda_i (u, w_i) = \sum_{i=1}^l \lambda_i m_0 = m_0.$$

Moreover, $(u, z) = m_0$ if and only if $\lambda_i > 0$ for every $i$ such that $(u, z_i) = (u, w_i)$. This happens if and only if $z \in Q$. This shows that $Q$ is a face of $\mathcal{P}_{(X,L)}$.

Let $\sigma \in \Sigma_X$ be an $n$-dimensional cone, and let $w$ be the corresponding vertex of $\mathcal{P}_{(X,L)}$. Then $\pi(w)$ is a vertex, say $v_1$, of $\pi_\mathbb{Q}(\mathcal{P}_{(X,L)})$ and hence $w$ lies on $R_1$. Since $R_1$ is also a face of $\mathcal{P}_{(X,L)}$, $w$ is a vertex of $R_1$. In each $R_i$, consider the vertex $w_i$ corresponding to the same cone of $\Sigma_Y$ as $w \in R_1$ (so we set $w_1 = w$). We have shown that $Q := \operatorname{Conv}(w_1, \ldots, w_l)$ is a face of $\mathcal{P}_{(X,L)}$, and $w = Q \cap R_1$.

Now call $\tau$ and $\eta$ the cones of $\Sigma_X$ corresponding respectively to $R_1$ and $Q$. It is $\sigma = \tau + \eta$, $\tau \subset \Delta_{\mathbb{Q}}$, and $\eta \cap \Delta_{\mathbb{Q}} = \{0\}$. This shows 3.2 $(ii)$ and thus $j$ induces a toric fibration. ∎

*Remark* 3.8 (Non reduced fibers). Let $f \colon X \to Y$ be a toric fibration. If $U_0 \subset Y$ is the open orbit of the torus action on $Y$, then $f^{-1}(U_0) \cong U_0 \times F$. If $X$ is smooth, then $Y$ and $F$ are smooth, and for every invariant, affine open subset $U$ of $Y$, we have $f^{-1}(U) \cong U \times F$. When $X$ is singular, $f$ can have non reduced fibers (even if $Y$ and $F$ are both smooth), and the property above is no longer true.

Let $L \in \operatorname{Pic} X$ be an ample line bundle, $p = V(\eta)$ a fixed point of $Y$, $V(\widetilde{\eta})$ the fiber of $f$ over $p$, and $Q$ the corresponding face of $\mathcal{P}_{(X,L)}$. Then the following are equivalent:

(i) the scheme-theoretical fiber of $f$ over $p$ is reduced;
(ii) $f^{-1}(U_\eta) \cong U_\eta \times F$;
(iii) $N = \Delta \oplus (\operatorname{Aff}(\widetilde{\eta}) \cap N)$;
(iv) $\pi_{|\operatorname{Aff}(Q) \cap M} \colon \operatorname{Aff}(Q) \cap M \to \Lambda$ is bijective.

The following lemma gives a local property of the polytope associated to a polarized toric fibration.

**Lemma 3.9.** *Let $(f \colon X \to Y, L)$ be a toric fibration polarized by a very ample line bundle $L \in \operatorname{Pic} X$. Suppose that $Q$ is a face of $\mathcal{P}_{(X,L)}$ such that:*

- *$Q$ has a vertex in the origin;*
- *$Q$ corresponds to an invariant fiber $F$ of $f$ such that the scheme-theoretical fiber of $f$ over $f(F)$ is reduced.*

*Then $Q = \mathcal{P}_{(F, L_{|F})}$, and there exists a basis $u_1, \ldots, u_k, u'_1, \ldots, u'_m$ of $M$ such that:*

- *$u_1, \ldots, u_k \in Q$;*
- *$u'_1, \ldots, u'_m \in \mathcal{P}_{(X,L)}$ and are a basis of $\ker \pi$;*
- *$u_h + u'_j \in \mathcal{P}_{(X,L)}$ for every $h = 1, \ldots, k$ and $j = 1, \ldots, m$.*

*Proof.* We have $\mathcal{P}_{(X,L)} = \mathcal{C}(R_1, \ldots, R_l, \pi)$ by Lemma 3.6. Let $p \in F$ be the fixed point corresponding to the origin, and let $R_1$ be the face of $\mathcal{P}_{(X,L)}$ corresponding to the invariant section of $f$ through $p$.

Set $\Gamma := \operatorname{Aff}(Q) \cap M$. Since the fiber of $f$ over $p$ is reduced, we know by Remark 3.8 that $\pi_{|\Gamma} \colon \Gamma \to \Lambda$ is an isomorphism and $M = \ker \pi \oplus \Gamma$.

It follows that $(\Gamma, Q)$ and $(\Lambda, \pi_{\mathbb{Q}}(P))$ are isomorphic lattice polytopes, hence $Q = \mathcal{P}_{(F, L_{|F})}$ by Lemma 3.6.

Let $w_1, \ldots, w_k$ be the vertices of $Q$ adjacent to the origin. Because $L_{|F}$ is very ample, Lemma 1.3 implies that the semigroup $\langle w_1, \ldots, w_k \rangle \cap \Gamma$ is generated by $Q \cap \Gamma$. Since the group generated by $\langle w_1, \ldots, w_k \rangle \cap \Gamma$ is the whole $\Gamma$, $Q$ must contain a basis $u_1, \ldots, u_k$ of $\Gamma$.

Observe that $\mathrm{Aff}(R_1) = \ker \pi_{\mathbb{Q}}$, and call $\tau$ the cone over $R_1$ in $\ker \pi_{\mathbb{Q}}$. By [Ewa96, Lemma V.3.5], the semigroup $\tau \cap \ker \pi$ has a unique minimal system of generators, which as above must contain a basis $u'_1, \ldots, u'_m$ of $\ker \pi$.

For every $i = 1, \ldots, l$ let $z_i := Q \cap R_i$ (note that $z_1, \ldots, z_l$ are all the vertices of $Q$, and $z_1 = 0$). We claim that $u'_j + z_i \in R_i$ for every $j = 1, \ldots, m$ and $i = 1, \ldots, l$. In fact, since the $R_i$'s are strictly combinatorially isomorphic, the cone over $R_i - z_i$ in $\ker \pi_{\mathbb{Q}}$ is $\tau$. Since $L_{|Y_i}$ is very ample, Lemma 1.3 implies that $\tau \cap \ker \pi$ has a system of generators contained in $R_i - z_i$. Such system of generators must contain the minimal one, and in particular it must contain $u'_1, \ldots, u'_m$.

Let $h \in \{1, \ldots, k\}$ and $j \in \{1, \ldots, m\}$. Since $u_h \in Q$, we can write $u_h = \sum_{i=1}^{l} \lambda_i z_i$, with $\lambda_i \geq 0$ and $\sum_{i=1}^{l} \lambda_i = 1$. Moreover for each $i$ we have $u'_j + z_i \in R_i$. Hence

$$u_h + u'_j = \left(\sum_{i=1}^{l} \lambda_i z_i\right) + u'_j = \sum_{i=1}^{l} \lambda_i \left(z_i + u'_j\right)$$

is a convex combination of points of $\mathcal{P}_{(X,L)}$. ∎

## 4. Geometric description

In this section we construct a natural fibration $\phi$ on $X \subset \mathbb{P}^N$ toric and $\mathbb{Q}$-factorial, covered by lines. The basic property of $\phi$ is that it contracts all families of lines covering $X$.

We are going to use the results of section §2 and the following result.

**Theorem 4.1** ([BCD07]). *Let $X \subset \mathbb{P}^N$ be a $\mathbb{Q}$-factorial toric variety, and let $V \subset G(1, N)$ be an irreducible subvariety parametrizing a covering family of lines in $X$.*

*Then there exists an elementary toric fibration $f \colon X \to Y$ such that every line of the family $V$ lies in some fiber of $f$.*

**Theorem 4.2.** *Let $X \subset \mathbb{P}^N$ be a $\mathbb{Q}$-factorial toric variety, and assume that $X$ is covered by lines. Then there exist (and are uniquely determined) a projective, $\mathbb{Q}$-factorial toric variety $Z$, and a surjective and equivariant morphism $\phi \colon X \to Z$, locally trivial in the Zariski topology, such that:*

(i) *the fiber $F$ is a product $J_1 \times \cdots \times J_r$;*

(ii) *each $J_i$ is a projective, $\mathbb{Q}$-factorial toric variety with $\rho_{J_i} = 1$. It is the projective join in $\mathbb{P}^N$ of pairwise disjoint invariant subvarieties;*

(iii) *each line in $X$ which intersects the open orbit of the torus, is contained in some fiber of $\phi$.*

*Proof.* Choose an irreducible family of lines which covers $X$, and let $f_1 \colon X \to Z_1$ be the elementary toric fibration given by Theorem 4.1. If every line contained in $X$ which intersects $U_0$ is contracted by $f_1$, we set $\phi := f_1$. Otherwise, there is a line $L \subset X$ which intersects $U_0$ and is not contracted by $f_1$. Moving $L$ with the action of the torus, we get a second irreducible family of lines covering $X$, not contracted by $f_1$. This produces a second elementary toric fibration $f_2 \colon X \to Z_2$.

Recall that every elementary toric fibration corresponds to a unique extremal ray of $\mathrm{NE}(X)$, and this cone is polyhedral, so there are a finite number of possibilities. Iterating the above procedure, we get $r$ distinct elementary toric fibrations $f_i \colon X \to Z_i$, such that for every $i = 1, \dots, r$ the general fiber $J_i$ is covered by lines, and such that every line in $X$ which intersects $U_0$ is contracted by some $f_i$.

By Corollary 2.8 we know that $\mathrm{def}(J_i) > 0$ and that each $J_i$ is the projective join of $\mathrm{def}(J_i) + 1$ pairwise disjoint invariant subvarieties.

Assume that $r \geq 2$ and let $\Delta_i \subset N$ be the sublattice corresponding to $f_i$, as in Lemma 3.2. We claim that $\Delta_i \cap \Delta_j = \{0\}$, for every $i, j \in \{1, \dots, r\}$. If so, it is easy to see that the conditions in Lemma 3.2 are satisfied for the sublattice $\Delta_1 \oplus \cdots \oplus \Delta_r \subseteq N$. This gives a toric fibration $\phi \colon X \to Z$, whose general fiber is $F = J_1 \times \cdots \times J_r$.

For simplicity we show that $\Delta_1 \cap \Delta_2 = \{0\}$. Let $v \in \Delta_1 \cap \Delta_2$ and let $\sigma \in \Sigma_X$ be an $n$-dimensional cone containing $v$. Let $\sigma = \tau_1 + \eta_1 = \tau_2 + \eta_2$ be the decompositions given by Lemma 3.2 $(ii)$, with the obvious notation. Then $V(\eta_1)$ (respectively, $V(\eta_2)$) is the invariant fiber of $f_1$ (respectively, $f_2$) through the point $V(\sigma)$. Since $f_1$ and $f_2$ are distinct elementary toric fibrations, we have $V(\eta_1) \cap V(\eta_2) = V(\sigma)$, equivalently $\eta_1 + \eta_2 = \sigma$. Write $\eta_1 = \eta_1 \cap \eta_2 + \gamma_1$ and $\eta_2 = \eta_1 \cap \eta_2 + \gamma_2$ with $\gamma_1, \gamma_2 \in \Sigma_X$, and $\gamma_1 \cap \gamma_2 = \{0\}$. Then we get $\gamma_1 = \tau_2$, $\gamma_2 = \tau_1$, and $\sigma = \tau_1 + \tau_2 + \eta_1 \cap \eta_2$. Hence $v \in \sigma \cap (\Delta_1)_{\mathbb{Q}} \cap (\Delta_2)_{\mathbb{Q}} = \tau_1 \cap \tau_2 = \{0\}$, and we are done.

We are left to show that $\phi$ is locally trivial in the Zariski topology. By remark 3.8, this is equivalent to showing that $\phi$ has no non reduced fibers. Suppose that there exists a fiber $F_1$ of multiplicity $m$. Then $(F_1)_{red} \cong F$ and there is a curve $C \subset F_1$ such that $mC$ is numerically equivalent to a line in $F$. This yields $m = 1$. ∎

## 5. Computing the dual defect

In this section we apply our results to compute the dual defect of $\mathbb{Q}$-factorial toric varieties in $\mathbb{P}^N$.

The exact formula for the defect will be derived by relating $\mathrm{def}(X)$ to the defect of the fiber of the fibration given by Theorem 4.2.

The following is a generalization of [LS87, Proposition 3.5].

**Proposition 5.1.** *Let $X \subset \mathbb{P}^N$ be a normal variety, and $f \colon X \to Y$ a morphism with connected fibers onto a projective variety $Y$. Assume that for some $q \in X^*_{reg}$ the contact locus of $H_q$ with $X$ lies in a fiber of $f$.*

*Then if $F$ is a general fiber of $f$, we have*

$$\mathrm{def}(X) \geq \mathrm{def}(F) - \dim Y.$$

*Moreover, in case of equality, there exists a morphism $f^* \colon X^*_{reg} \to Y$ such that for a general point $y \in Y$, the closure of $(f^*)^{-1}(y)$ in $X^*$ is $(f^{-1}(y))^*$.*

*Proof.* Consider diagram (1) in section 1, and let $\mathcal{J}_X \xrightarrow{g} Z \xrightarrow{h} X^*$ be the Stein factorization of $\pi_2$. Recall that $\mathcal{J}_{X^*}^0 = \pi_2^{-1}(X^*_{reg})$. Because every fiber of $\pi_2$ over $X^*_{reg}$ is a $\mathbb{P}^{\operatorname{def}(X)}$, $h$ is an isomorphism over $X^*_{reg}$, and every fiber of $g$ over $h^{-1}(X^*_{reg})$ is a $\mathbb{P}^{\operatorname{def}(X)}$. By hypothesis $f \circ \pi_1$ contracts one of these fibers, so the same must hold for all.

By [Deb01, Lemma 1.15] there exists a morphism $f^*\colon X^*_{reg} \to Y$ such that the following diagram commutes:

$$\begin{array}{ccc}
 & \mathcal{J}_{X^*}^0 & \\
{\scriptstyle (\pi_1)_{|\mathcal{J}_{X^*}^0}} \swarrow & & \searrow {\scriptstyle (\pi_2)_{|\mathcal{J}_{X^*}^0}} \\
X & & X^*_{reg} \\
{\scriptstyle f} \searrow & & \swarrow {\scriptstyle f^*} \\
 & Y &
\end{array}$$

Let $y_0 \in Y$ be a general point and set $F_{y_0} := f^{-1}(y_0)$. The fiber $F_{y_0} \subset \mathbb{P}^N$ is an integral variety, and $(f^*)^{-1}(y_0)$ is smooth and connected. We claim that:
- $(f^*)^{-1}(y_0) \subseteq (F_{y_0})^*$;
- if $q \in (f^*)^{-1}(y_0)$ is a general point, the contact locus of $H_q$ with $F_{y_0}$ contains the contact locus of $H_q$ with $X$.

In fact, let $q \in (f^*)^{-1}(y_0)$ be a general point and let $H_q$ be the corresponding hyperplane. Let $M_0 := \pi_1(\pi_2^{-1}(q)) \cong \mathbb{P}^{\operatorname{def}(X)}$ be the contact locus of $H_q$ with $X$. By the commutativity of the diagram above, $M_0 \subseteq F_{y_0}$. By the generality of $y_0$ and $q$, $M_0 \cap X_{reg} \cap (F_{y_0})_{reg} \neq \emptyset$. For any $x \in M_0 \cap X_{reg} \cap (F_{y_0})_{reg}$ we have $H_q \supseteq \mathbb{T}_x X \supseteq \mathbb{T}_x F_{y_0}$ and thus $q \in (F_{y_0})^*$ and the contact locus of $H_q$ with $F_{y_0}$ contains $M_0$.

Since $f^*$ is dominant it must be:
$$\dim F^*_{y_0} \geq \dim X^* - \dim Y,$$
namely
$$\operatorname{def}(X) \geq \operatorname{def}(F_{y_0}) - \dim Y.$$

Suppose that $\operatorname{def}(X) = \operatorname{def}(F_{y_0}) - \dim Y$. Then $\dim(F_{y_0})^* = \dim X^* - \dim Y$. We have $(f^*)^{-1}(y_0) \subseteq (F_{y_0})^*$ and $\dim(f^*)^{-1}(y_0) = \dim(F_{y_0})^*$, so the closure of $(f^*)^{-1}(y_0)$ in $X^*$ is the dual variety $(F_{y_0})^*$. ∎

Let $X \subset \mathbb{P}^N$ be an embedding defined by the complete linear system $|\mathcal{O}_X(1)|$, and let $p_0 \in X$ be a general point.

Choose $n$ sections in $H^0(X, \mathcal{O}_X(1))$ which give local coordinates $x_1, \ldots, x_n$ in $p_0$. Let $g \in H^0(X, \mathcal{O}_X(1))$ be a section vanishing in $p_0$, and view locally $g$ as a function $g(x_1, \ldots, x_n)$. We denote by $\operatorname{Hes}_{p_0}(g)$ the Hessian matrix of $g$ in $p_0$. Then we have (see [Kat73] and [GKZ94, §1.5]):

(4) $$\operatorname{def}(X) = \dim X - \max_g (\operatorname{rk} \operatorname{Hes}_{p_0}(g)),$$

where $g$ varies among global sections of $\mathcal{O}_X(1)$ vanishing in $p_0$.

The local description given in Lemma 3.9 turns out to be very useful in order to apply (4).

**Theorem 5.2.** *Let $X \subset \mathbb{P}^N$ be a $\mathbb{Q}$-factorial toric variety covered by lines. Let $\phi\colon X \to Z$ be the toric fibration given by Theorem 4.2, and $F \subset \mathbb{P}^N$ a fiber. Then*
$$\operatorname{def}(X) = \max(0, \operatorname{def}(F) - \dim Z).$$
*Moreover, if $\operatorname{def}(X) = \operatorname{def}(F) - \dim Z$, then there exists a morphism $\phi^*\colon X^*_{reg} \to Z$ such that for a general point $z \in Z$, the closure of $(\phi^*)^{-1}(z)$ in $X^*$ is $(\phi^{-1}(z))^*$.*

Notice that in particular, if $X$ has positive dual defect, then $\dim Z < \frac{1}{2}\dim X$.

*Proof.* First of all observe that we can assume that $X$ is embedded in $\mathbb{P}^N$ by the complete linear system $|\mathcal{O}_X(1)|$. In fact, this does not change $\operatorname{def}(X)$ nor $\operatorname{def}(F)$ by Remark 1.1.

Set $k := \dim F$ and $m := \dim Z$. Let $p_0$ be the distinguished point of the torus $T$ (i.e. the point $(1,\ldots,1)$), and $F_0$ the fiber of $\phi$ through $p_0$. Notice that since $\phi$ is locally trivial in the Zariski topology, for every fiber $F$ of $\phi$ there is an isomorphism $\xi\colon F \to F_0$ such that $\xi^*(\mathcal{O}_{F_0}(1)) = \mathcal{O}_F(1)$, therefore $\operatorname{def}(F) = \operatorname{def}(F_0)$.

The restriction $H^0(X, \mathcal{O}_X(1)) \to H^0(F_0, \mathcal{O}_{F_0}(1))$ is surjective, hence $F_0$ is embedded by the complete linear system $|\mathcal{O}_{F_0}(1)|$.

Up to translation, we can assume that the polytope $\mathcal{P}_X$ has a vertex in the origin. Let $Q$ be the face of $\mathcal{P}_X$ containing the origin and corresponding to an invariant fiber of $\phi$. By Theorem 4.2, all fibers of $\phi$ are reduced. In particular Lemma 3.9 applies to $\mathcal{P}_X$ and $Q$. This says that $Q = \mathcal{P}_{F_0}$.

Recall that $H^0(X, \mathcal{O}_X(1)) = \oplus_{u \in \mathcal{P}_X \cap M}\mathbb{C}\chi^u$, where $\chi^u$ is the rational function on $X$ corresponding to $u$. Let $u_1,\ldots,u_k,u'_1,\ldots,u'_m$ be the basis of $M$ given by Lemma 3.9, and set
$$x_i := \chi^{u_i} - 1, \qquad x'_j := \chi^{u'_j} - 1$$
for all $i = 1,\ldots,k$ and $j = 1,\ldots,m$. Then $x_1,\ldots,x_k,x'_1,\ldots,x'_m \in H^0(X, \mathcal{O}_X(1))$ are such that:
- $x_1,\ldots,x_k,x'_1,\ldots,x'_m$ are local coordinates in $p_0$;
- locally $F_0$ is given by $x'_1 = \cdots = x'_m = 0$;
- $x_h x'_j \in H^0(X, \mathcal{O}_X(1))$ for any $h = 1,\ldots,k$ and $j = 1,\ldots,m$.

Since $H^0(F_0, \mathcal{O}_{F_0}(1)) = \oplus_{u \in Q \cap M}\mathbb{C}\chi^u$, any $g \in H^0(F_0, \mathcal{O}_{F_0}(1))$ can be written as $g = \sum_j \lambda_j \chi^{w_j}$, with $w_j \in Q$. By our choice of $u_1,\ldots,u_k$, we can view $g$ as a function of $x_1,\ldots,x_k$. By (4), we can choose such a $g$ with $r_0 := \operatorname{rk}\operatorname{Hes}_{p_0}(g) = k - \operatorname{def}(F)$.

Choose a symmetric minor $A$ of $\operatorname{Hes}_{p_0}(g)$ of order $r_0$ with non zero determinant, and let $j_1,\ldots,j_{\operatorname{def}(F)}$ be the indices of the rows not appearing in $A$.

Set $l := \min(m, \mathrm{def}(F))$ and consider

$$g' := g + \sum_{i=1}^{l} x_{j_i} x'_i.$$

Then $g'$ is a global section of $\mathcal{O}_X(1)$, vanishing in $p_0$. The following figure shows $\mathrm{Hes}_{p_0}(g')$ in the case $j_1 = r+1, \ldots, j_{\mathrm{def}(F)} = k$.

$$k \left\{ \begin{bmatrix} A & * & 0 \\ & & I_l \\ * & * & \\ & & 0 \\ \hline 0 & I_l & 0 & 0 \end{bmatrix} \right. \quad \begin{bmatrix} A & * & 0 \\ * & * & I_l \\ & I_l & \\ 0 & & 0 \end{bmatrix} \left. \begin{matrix} \} k \\ \\ \} m \end{matrix} \right.$$

$$\text{Case } l = m \qquad \qquad \text{Case } l = \mathrm{def}(F)$$

Suppose that $\mathrm{def}(F) \geq \dim Z = m$, so that $l = m$. Consider the symmetric minor of $\mathrm{Hes}_{p_0}(g')$ containing the $r_0$ rows of $A$, the rows $j_1, \ldots, j_m$, and the last $m$ rows. This minor has non zero determinant, so $\mathrm{rk}\,\mathrm{Hes}_{p_0}(g') \geq r_0 + 2m$. By (4) this gives

$$\mathrm{def}(X) \leq n - r_0 - 2m = \mathrm{def}(F) - \dim Z.$$

On the other hand, we know that $\phi$ contracts to a point every linear subspace contained in $X$ passing through $p_0$. Therefore Proposition 5.1 implies that $\mathrm{def}(X) = \mathrm{def}(F) - \dim Z$.

Suppose now that $\mathrm{def}(F) < \dim Z$, so that $l = \mathrm{def}(F)$. Consider the symmetric minor of $\mathrm{Hes}_{p_0}(g')$ containing the first $k + \mathrm{def}(F)$ rows. This minor has non zero determinant, so $\mathrm{rk}\,\mathrm{Hes}_{p_0}(g') \geq k + \mathrm{def}(F)$. By (4) this gives $\mathrm{def}(X) \leq n - k - \mathrm{def}(F) = \dim Z - \mathrm{def}(F)$, which is not enough.

We proceed by contradiction and assume that $\mathrm{def}(X) > 0$. Then by (4) we have $\mathrm{rk}\,\mathrm{Hes}_{p_0}(g') < n$. The kernel of $\mathrm{Hes}_{p_0}(g')$ is positive dimensional, and is contained in $x_1 = \cdots = x_k = x'_1 = \cdots = x'_{\mathrm{def}(F)} = 0$. Observe that we can also assume that $g'$ vanishes with order at least two in $p_0$. This means that $g'$ gives a hyperplane $H_{q_0} \subset \mathbb{P}^N$ such that $H_{q_0} \cap X$ is singular at $p_0$ (so $q_0 \in X^*$), and locally its singular locus *is not* contained in $F_0$.

Theorem 4.2 implies that every line in $X$ through $p_0$ is contained in $F_0$, so $\mathrm{Sing}(H_{q_0} \cap X)$ contains no lines through $p_0$. Consider diagram (1) of section 1. For a general $q \in X^*$, $\pi_1(\pi_2^{-1}(q))$ is a linear $\mathbb{P}^{\mathrm{def}(X)}$ in $\mathbb{P}^N$, hence for *every* $q \in X^*$, the contact locus $\pi_1(\pi_2^{-1}(q))$ must be covered by lines.

Observe that $\pi_1(\pi_2^{-1}(q_0)) \cap X_{reg} = \mathrm{Sing}(H_{q_0} \cap X) \cap X_{reg}$, hence every irreducible component of $\pi_1(\pi_2^{-1}(q_0))$ containing $p_0$ must be contained in $\mathrm{Sing}(H_{q_0} \cap X)$. Therefore $\pi_1(\pi_2^{-1}(q_0))$ contains no line passing through $p_0$, which is a contradiction.

The last part of the statement follows from Proposition 5.1. ∎

*Remark* 5.3. The proof of Theorem 5.2 shows more generally the following. Let $X \subset \mathbb{P}^N$ be a Q-factorial toric variety endowed with a toric fibration $f\colon X \to Y$. Suppose that there is at least one fixed point $p \in Y$ such that the scheme-theoretical fiber of $f$ over $p$ is reduced. Then

$$\operatorname{def}(X) \leq |\operatorname{def}(F) - \dim Y|,$$

where $F$ is a general fiber of $f$.

**Corollary 5.4.** *Let $X \subset \mathbb{P}^N$ be a $\mathbb{Q}$-factorial toric variety covered by lines, and let $\phi$ be the fibration given by Theorem 4.2.*

*The dual defects of $F = J_1 \times \cdots \times J_r$ and $X$ are given by the following formulas:*

$$\operatorname{def}(F) = \max(0, \dim J_1 + \operatorname{def}(J_1) - \dim F, \ldots, \dim J_r + \operatorname{def}(J_r) - \dim F),$$
$$\operatorname{def}(X) = \max(0, \dim J_1 + \operatorname{def}(J_1) - \dim X, \ldots, \dim J_r + \operatorname{def}(J_r) - \dim X).$$

*Proof.* The dual defect of $F = J_1 \times \cdots \times J_r$ is given by [WZ94, Theorem 3.1]. Moreover Theorem 5.2 gives $\operatorname{def}(X) = \max(0, \operatorname{def}(F) - \dim Z)$ and hence the statement. ∎

In [DR06] it is shown that a smooth toric $X \subset \mathbb{P}^N$ has positive dual defect if and only if $X$ is a $\mathbb{P}^k$-bundle with $k > \frac{1}{2}\dim X$. The following is the generalization to the $\mathbb{Q}$-factorial case.

**Corollary 5.5.** *Let $X \subset \mathbb{P}^N$ be a toric and $\mathbb{Q}$-factorial variety. Then $X$ has positive dual defect if and only if there exists an elementary toric fibration $\psi\colon X \to Y$ whose fiber is a projective join $J$ with $\dim J + \operatorname{def}(J) > \dim X$.*

*Proof.* Assume that $\operatorname{def}(X) > 0$. Then $X$ is covered by lines and Theorem 4.2 implies the existence of a toric fibration $\phi\colon X \to Z$ with fiber $F = J_1 \times \cdots \times J_r$, $r \geq 1$. By Corollary 5.4, $\operatorname{def}(X) > 0$ implies that $\dim J_{j_0} + \operatorname{def}(J_{j_0}) > \dim X$ for some $j_0 \in \{1, \ldots, r\}$. The numerical class of a line in $J_{j_0}$ is extremal in $\operatorname{NE}(X)$, and determines an elementary toric fibration $\psi\colon X \to Y$ with fiber $J_{j_0}$.

Conversely, suppose that there exists $\psi$ as in the statement. Then $X$ is covered by lines, and again Theorem 4.2 implies the existence of a toric fibration $\phi\colon X \to Z$ with fiber $F = J_1 \times \cdots \times J_r$. Since $\phi$ contracts all lines in $X$ passing through a general point, $J$ must be one of $J_1, \ldots, J_r$. Then Corollary 5.4 implies that $\operatorname{def}(X) > 0$. ∎

*Remark* 5.6. In the situation of Corollary 5.5, one can always choose a suitable $\psi$ such that $\operatorname{def}(J) = \operatorname{def}(X) + \dim Y$. For related results when $X \subset \mathbb{P}^N$ is smooth and possibly non toric, see [BFS92, Theorem (1.2)].

*Example* 5.7 (Low dimensions). In the following list we describe all toric and $\mathbb{Q}$-factorial varieties $X \subset \mathbb{P}^N$, covered by lines, of dimension $n \leq 3$. In the table $Z$ is the basis of the fibration $\phi$ given by Theorem 4.2.

| $n$ | $\rho_X$ | $X$ | $Z$ | $\mathrm{def}(X)$ |
|---|---|---|---|---|
| 1 | 1 | $\mathbb{P}^1$ | $pt$ | 1 |
| 2 | 1 | $\mathbb{P}^2$ | $pt$ | 2 |
| 2 | 1 | $J$ cone over a non linear rational normal curve | $pt$ | 1 |
| 2 | 2 | $\mathbb{P}^1 \times \mathbb{P}^1 \subset \mathbb{P}^3 \subseteq \mathbb{P}^N$ | $pt$ | 0 |
| 2 | 2 | $\mathbb{P}^1$-bundle over $\mathbb{P}^1$ | $\mathbb{P}^1$ | 0 |
| 3 | 1 | $\mathbb{P}^3$ | $pt$ | 3 |
| 3 | 1 | cone over a non linear rational normal curve with vertex a line | $pt$ | 2 |
| 3 | 1 | cone over a $\mathbb{Q}$-factorial toric surface not covered by lines | $pt$ | 1 |
| 3 | 1 | join of two disjoint non linear rational normal curves | $pt$ | 1 |
| 3 | 2 | $\mathbb{P}^2 \times \mathbb{P}^1 \subset \mathbb{P}^5 \subseteq \mathbb{P}^N$ | $pt$ | 1 |
| 3 | 2 | $X = J \times \mathbb{P}^1$, $J$ as above | $pt$ | 0 |
| 3 | 3 | $\mathbb{P}^1 \times \mathbb{P}^1 \times \mathbb{P}^1 \subset \mathbb{P}^7 \subseteq \mathbb{P}^N$ | $pt$ | 0 |
| 3 | 2 | $\mathbb{P}^2$-bundle over $\mathbb{P}^1$ | $\mathbb{P}^1$ | 1 |
| 3 | 2 | fibration in $J$ over $\mathbb{P}^1$, $J$ as above | $\mathbb{P}^1$ | 0 |
| 3 | 3 | $(\mathbb{P}^1 \times \mathbb{P}^1)$-bundle over $\mathbb{P}^1$ | $\mathbb{P}^1$ | 0 |
| 3 | $1 + \rho_S$ | $\mathbb{P}^1$-bundle over a $\mathbb{Q}$-factorial toric surface $S$ | $S$ | 0 |

## 6. Combinatorial characterization

The dual defect of a toric variety plays an important role in the study of the $A$-discriminant, associated to a subset $A \subset \mathbb{Z}^n$. We refer to [GKZ94] for details.

A finite subset $A \subset \mathbb{Z}^n$ is viewed as a collection of exponents of Laurent monomials in $n$ variables, so it defines a map $f_A : (\mathbb{C}^*)^n \to \mathbb{P}^{|A|-1}$. We denote the closure of the image by $X_A$. If $X_A$ is normal, then it is a toric variety. Otherwise it can be considered as a toric variety in a generalized sense, see [GKZ94, §5.1].

When the dual variety $(X_A)^*$ is a hypersurface, there exists an irreducible homogeneous polynomial $F_A$, with integral coefficients, defining $(X_A)^*$.

The $A$-discriminant $D_A$ is defined as:

$$D_A = \begin{cases} F_A & \text{if } \mathrm{def}(X_A) = 0 \\ 1 & \text{if } \mathrm{def}(X_A) > 0. \end{cases}$$

In [DS02] the case where $X_A$ has codimension two in $\mathbb{P}^{|A|-1}$ is studied. The authors give an explicit formula for $D_A$ and a combinatorial condition for $D_A = 1$. Recently in [CC07] the combinatorial characterization of defectivity has been generalized to the case $\mathrm{codim} X_A \leq 4$.

As a corollary of our results we obtain a characterization in any codimension, under suitable hypotheses on the polytope $\mathrm{Conv}(A)$. The following shows that, when $D_A = 1$, the polytope $\mathrm{Conv}(A)$ has the structure of a $\pi$-twisted Cayley sum.

**Theorem 6.1.** *Let $A \subset \mathbb{Z}^n$ be a finite subset such that $\mathrm{Conv}(A) \subset \mathbb{Q}^n$ is a simple polytope of dimension $n$.*

*For every vertex $v$ of $\mathrm{Conv}(A)$, let $P_v$ be the translated polytope $\mathrm{Conv}(A) - v$, and let $\eta_v$ be the cone over $P_v$ in $\mathbb{Q}^n$. Assume that the semigroup $\eta_v \cap \mathbb{Z}^n$ is generated by $A - v$.*

*Then* $D_A = 1$ *if and only if there exists a surjection of lattices* $\pi\colon \mathbb{Z}^n \to \Lambda$ *such that:*

(i) $S := \operatorname{Conv}(\pi(A)) \subset \Lambda_{\mathbb{Q}}$ *is a lattice simplex* $\operatorname{Conv}(v_0, \ldots, v_{\dim S})$ *with* $\dim S + \operatorname{def}(S) > n$;

(ii) $A \subset \operatorname{Conv}(A_0, \ldots, A_{\dim S})$, *where* $A_i := A \cap \pi^{-1}(v_i)$;

(iii) *the polytopes* $\operatorname{Conv}(A_i) \subset \mathbb{Q}^n$ *are strictly combinatorially isomorphic for* $i = 0, \ldots, \dim S$.

*Proof.* The lattice polytope $\operatorname{Conv}(A) \subset \mathbb{Q}^n$ determines a pair $(X, L)$ where $X$ is a toric variety of dimension $\dim \operatorname{Conv}(A) = n$ and $L \in \operatorname{Pic} X$ is ample. Notice that $X$ is $\mathbb{Q}$-factorial because $\operatorname{Conv}(A)$ is simple. Moreover the hypotheses on $A$ say that the sublinear system $|V|$ of $|L|$ given by

$$V := \bigoplus_{u \in A} \mathbb{C}\chi^u \subseteq H^0(X, L)$$

determines an embedding $\varphi_{|V|}\colon X \hookrightarrow \mathbb{P}^{|A|-1}$. The restriction of $\varphi_{|V|}$ to the open subset $(\mathbb{C}^*)^n \subset X$ coincides with $f_A$, so $X_A$ is the image of $\varphi_{|V|}$. In particular it is isomorphic to $X$, and $\mathcal{P}_{X_A} = \operatorname{Conv}(A)$.

Assume that $D_A = 1$. By definition $\operatorname{def}(X_A) > 0$ and thus Corollary 5.5 implies the existence of an elementary toric fibration $\psi\colon X_A \to Y$ with fiber a projective join $J$ such that $\dim J + \operatorname{def}(J) > n$. By Lemma 3.6, $\operatorname{Conv}(A)$ is a $\pi$-twisted Cayley sum $\mathcal{C}(R_1, \ldots, R_l, \pi)$, where $\pi$ is dual to the map inducing the fibration $\psi$. The polytope $S = \pi_{\mathbb{Q}}(\operatorname{Conv}(A))$ is the polytope associated to $J \subset \mathbb{P}^{|A|-1}$, so it a simplex of dimension $\dim J$. Corollary 2.9 gives $\operatorname{def}(S) = \operatorname{def}(J)$, so we have (i). Moreover $\pi_{\mathbb{Q}}(R_1), \ldots, \pi_{\mathbb{Q}}(R_l)$ are the vertices of $S$, and $\operatorname{Conv}(A) = \operatorname{Conv}(R_1, \ldots, R_l)$, so we have (ii) and (iii).

Conversely, assume that there exists $\pi\colon \mathbb{Z}^n \to \Lambda$ as in the statement, and set $R_i := \operatorname{Conv}(A_i)$ for $i = 0, \ldots, \dim S$. Then $\operatorname{Conv}(A)$ is the $\pi$-twisted Cayley sum of $R_0, \ldots, R_{\dim S}$. Hence Lemma 3.6 implies the existence of a toric fibration $\psi\colon X_A \to Y$ with general fiber $J$ such that $\mathcal{P}_J = S$. We have $\rho_J = 1$ because $S$ is a simplex. Hence Corollary 2.9 yields $\operatorname{def}(J) = \operatorname{def}(S)$. Then we get $\dim J + \operatorname{def}(J) > n$ and Corollary 5.5 implies that $\operatorname{def}(X_A) > 0$. ∎

When $\operatorname{Conv}(A)$ is a simplex, the conditions in Theorem 6.1 become particularly simple. Geometrically this corresponds to the case of Picard number 1.

**Corollary 6.2.** *Let* $A \subset \mathbb{Z}^n$ *be a finite subset such that* $\operatorname{Conv}(A) \subset \mathbb{Q}^n$ *is an $n$-dimensional simplex.*

*For every vertex $v$ of* $\operatorname{Conv}(A)$, *let $P_v$ be the translated polytope* $\operatorname{Conv}(A) - v$, *and let $\eta_v$ be the cone over $P_v$ in $\mathbb{Q}^n$. Assume that the semigroup $\eta_v \cap \mathbb{Z}^n$ is generated by $A - v$.*

*Then the following are equivalent:*

(i) $D_A = 1$;

(ii) $\operatorname{def}(\operatorname{Conv}(A)) > 0$;

(iii) there exist two disjoint faces $Q_0$, $Q_1$ of $\mathrm{Conv}(A)$ such that $\mathrm{Conv}(A) \cap \mathbb{Z}^n = (Q_0 \cap \mathbb{Z}^n) \cup (Q_1 \cap \mathbb{Z}^n)$.

*Proof.* Suppose that $R_0$ and $R_1$ are two disjoint faces of $\mathrm{Conv}(A)$ having the same dimension, and such that $\mathrm{Aff}(R_0)$ is parallel to $\mathrm{Aff}(R_1)$. Since $\mathrm{Conv}(A)$ is a simplex, the only possibility is that $R_i$ are points. This means that the only possibility for a surjection $\pi$ as in Theorem 6.1 is $\pi = \mathrm{Id}_{\mathbb{Z}^n}$. Then Theorem 6.1 says that $D_A = 1$ if and only if $\mathrm{Conv}(A)$ has positive lattice defect. By the definition of lattice defect (see 2.5), this is also equivalent to (iii). ∎

The following provides examples of polytopes $\mathrm{Conv}(A)$ with $D_A = 1$ and $D_A \neq 1$.

*Example* 6.3. Consider the following lattice simplices:

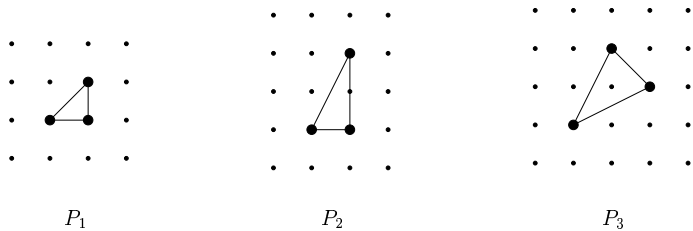

$P_1$        $P_2$        $P_3$

It is easy to see that $\mathrm{def}(P_1) = 2$, $\mathrm{def}(P_2) = 1$, and $\mathrm{def}(P_3) = 0$. Observe that $P_1$ is a standard simplex.

Suppose that $\mathrm{Conv}(A)$ is as in Theorem 6.1, with $D_A = 1$, and it is not a simplex. Observe that $\mathrm{def}(S) \leq \dim S$, hence $\mathrm{rk}\,\Lambda = \dim S > \frac{1}{2}n$. In particular, $\Lambda$ can never be a trivial lattice, and $n \geq 3$. If $n = 3$, the corresponding toric variety $X_A$ is always smooth. For examples of this situation we refer the reader to [DR06, Examples 2 and 4].

*Example* 6.4. Let $P$ be the convex hull in $\mathbb{Q}^4$ of the following points:

$$(0,0,0,0),\ (1,0,0,0),\ (0,0,1,0),\ (0,2,1,0),$$
$$(0,0,0,1),\ (1,0,0,1),\ (0,0,1,1),\ (0,2,1,1).$$

Observe that $|P \cap \mathbb{Z}^4| = 10$ and that geometrically $P$ corresponds to the embedding $Q \times \mathbb{P}^1 \subset \mathbb{P}^9$, where $Q \subset \mathbb{P}^4$ is a quadric with singular locus a line. Set $A := P \cap \mathbb{Z}^4$. It is not difficult to see that $A$ satisfies the hypotheses of Theorem 6.1.

Consider the projection $\pi \colon \mathbb{Z}^4 \to \mathbb{Z}^3$ on the first three coordinates. Then $\pi_{\mathbb{Q}}(P)$ is a 3-dimensional simplex with lattice defect 2, so condition (i) of Theorem 6.1 is satisfied. It is easy to see that also (ii) and (iii) are fulfilled, hence $D_A = 1$.

*Example* 6.5. Consider the polytope $P$ described in example 3.7 and $A := P \cap M$. The hypotheses of Theorem 6.1 are fulfilled. Since $P$ is not a simplex, and it has only two parallel faces, the only possibility for a surjection $\pi$ as in Theorem 6.1 is the projection defined in example 3.7. For this map, conditions (ii) and (iii) are satisfied. However $\pi_{\mathbb{Q}}(P)$ has dimension 1 and lattice defect 1, hence (i) is not satisfied. Thus Theorem 6.1 implies that $D_A \neq 1$.

Università di Pisa - Dipartimento di Matematica - Largo B. Pontecorvo, 5 I-56127 Pisa, Italy
  *E-mail address*: casagrande@dm.unipi.it

K.T.H. Matematik - S-10044 Stockholm, Sweden
  *E-mail address*: dirocco@math.kth.se